\documentclass[11pt]{article}
\usepackage{latexsym,amsmath,amssymb,amscd}

\usepackage{amsmath}

\makeatletter \@addtoreset{figure}{section}
\def\thefigure{\thesection.\@arabic\c@figure}
\def\fps@figure{h,t}
\@addtoreset{table}{bsection}

\def\thetable{\thesection.\@arabic\c@table}
\def\fps@table{h, t}
\@addtoreset{equation}{section}

\makeatother

\newenvironment{proof}[1][Proof]{\textbf{#1.} }{\ \rule{0.5em}{0.5em}}

\newcommand{\R}{{\mathbb  R}}

\newcommand{\ds}{\displaystyle}
\newcommand{\ol}{\overline}
\numberwithin{equation}{section}
\newtheorem{thm}{\bf Theorem}[section]
\newtheorem{lem}[thm]{\bf Lemma}

\newcommand{\p}{\prime}
\newcommand{\de}{\delta}
\newcommand{\va}{\varphi}

\newcommand{\pa}{\partial}

\newcommand{\al}{\alpha}

\begin{document}

\newtheorem{theorem}{Theorem}[section]
\newtheorem{definition}[theorem]{Definition}
\newtheorem{lemma}[theorem]{Lemma}
\newtheorem{remark}[theorem]{Remark}
\newtheorem{proposition}[theorem]{Proposition}
\newtheorem{corollary}[theorem]{Corollary}
\newtheorem{example}[theorem]{Example}

\newcommand{\bfi}{\bfseries\itshape}

\newsavebox{\savepar}
\newenvironment{
boxit}{\begin{lrbox}{\savepar}
\begin{minipage}[b]{15.5cm}}{\end{minipage}\end{lrbox}
\fbox{\usebox{\savepar}}}

\title{{\bf Equivalence of energy methods in stability theory}}
\author{Petre Birtea, Mircea Puta}
\date{}
\maketitle

\begin{abstract}
We will prove the equivalence of three methods, the so called
energy methods, in order to establish the stability of an
equilibrium point for a dynamical system. We will illustrate by
examples that this result simplifies enormously the amount of
computations especially when the stability cannot be decided with
one of the three methods.
\end{abstract}



\section{Introduction.}

Let $M$ be a smooth manifold and \begin{equation}\label{11}\dot x
= f(x)\end{equation} be a dynamical system on $M$ given by the
vector field $f \in \mathfrak{X}\left( M\right)$ and suppose $x_e
\in M$ is an equilibrium state for (\ref{11}), i.e. $f(x_e) = 0$.
The problem of nonlinear stability of equilibrium states is a very
old one and the most know and remarkable results were obtained by
Lyapunov \cite{4}. They are based on finding what is called a
Lyapunov function $V \in C^1 (M, \mathbb{R})$ such that:
\begin{itemize}
  \item [(i)] $V(x_e) = 0$
  \item [(ii)] $V(x) >0$, for $x \not = x_e$
  \item [(iii)] $\dot{V} \leq 0$, where $\dot{V}$ is the derivative of $V$ along
  the trajectories of (\ref{11}).
\end{itemize}

In practice it is sometimes very difficult to find such a
function. In many situations one can use constants of  motion as
Lyapunov functions, i.e. functions $V:M\rightarrow \mathbb{R}$
such that $\dot{V}=0$. This was extensively used in the context of
Hamilton-Poisson systems where the Hamiltonian and the Casimirs of
the Poisson structure are constants of motion. The methods for
studying stability using constants of motion are the so called
energy methods. The most general results using this methods for
establishing stability can be found in \cite{pat} and \cite{vic}.
Since in the present paper we are discussing local nonlinear
stability we can replace, by considering a coordinate chart around
the equilibrium $x_{e}$, the manifold $M$ with $\mathbb{R}^{n}$,
where $n$ is the dimension of $M$.

In 1965 Arnold \cite{1} gives the following criteria for
determining nonlinear stability for an equilibrium point of
(\ref{11}).

\begin{theorem}\label{tA}
{\rm (The Arnold method \cite{1})} Let $C_1, \dots , C_k \in C^2
(\mathbb{R}^n, \mathbb{R})$ be constants of motion for the
equation (\ref{11}) and $F_i \in C^2 (\mathbb{R}^n \times
\mathbb{R}^{k-1}, \mathbb{R})$ be the smooth function given by:
\begin{equation*}
    F_i (x, \lambda_1, \dots, \widehat{\lambda}_i, \dots, \lambda_k)
    :\stackrel{def}{=} C_i (x)-\lambda_1 C_1(x)- \dots -\widehat{\lambda_i C_i
    (x)}- \dots -\lambda_k C_k(x)
\end{equation*} where $ \widehat{g}$ means that the term $g$ is
omitted. If there exist constants $\lambda_1^*, \dots ,$
$\widehat{\lambda_i^*}, \dots, \lambda_k^*$ in $\mathbb{R}$ such
that
\begin{itemize}
  \item [(i)] $\nabla_x F_i (x_e, \lambda_1^*, \dots, \widehat{\lambda_i^*},
  \dots, \lambda_k^*)=0$
  \item [(ii)] $\nabla_{xx}^2 F_i (x_e, \lambda_1^*, \dots,
  \widehat{\lambda_i^*}, \dots, \lambda_k^*)_{|_{W\times W}}$ is positive or
  negative definite, where $$W := \bigcap\limits_{\tiny \begin{array}{l}
  j=1\\j \not = i\end{array}}^k \ker d C_j
  (x_e),$$
\end{itemize}
then $x_e$ is nonlinear stable.
\end{theorem}

Later, in 1985, Holm, Marsden, Ratiu and Weinstein \cite{3} give
another method for establishing stability of an Hamilton-Poisson
system, the so called Energy-Casimir method.

\begin{theorem}\label{tB}{\rm (The Energy-Casimir method \cite{3})} Let
$C_1, C_2, \dots, C_k \in C^2 (\mathbb{R}^n, \R)$ be constants of
motion for the equation (\ref{11}). If there exist $\va_1, \dots,
\widehat{\va_i}, \dots, \va_k\in C^2 (\R, \R)$ such that:
\begin{itemize}
  \item [(i)] $\nabla_x (C_i + \va_1 (C_1) + \dots +
  \widehat{\va_i(C_i)} +\dots + \va_k (C_k))(x_e)=0$
  \item [(ii)] $\nabla_{xx}^2 (C_i +\va_1 (C_1) + \dots + \widehat{\va_i
  (C_i)}+\dots + \va_k (C_k))(x_e)$ is positive or negative
  definite,
\end{itemize}
then $x_e$ is nonlinear stable.
\end{theorem}

The above result has also an infinite dimensional analogue for
Hamilton-Poison systems on Banach spaces, see \cite{3}.

Studying the stability of relative equilibria, in 1998, Ortega and
Ratiu \cite{5} obtain, as a corollary of their results about
stability of relative equilibria, the following theorem.

\begin{theorem}\label{tC} {\rm (The Ortega-Ratiu method \cite{5})} Let $C_1,
\dots, C_k \in C^2 (\R^n, \R)$ be constants of motion for the
equation (\ref{11}). If there exist $\va_1, \dots,
\widehat{\va_i}, \dots, \va_k \in C^2 (\R, \R)$ such that:
\begin{itemize}
  \item [(i)] $\nabla_x (C_i + \va_1 (C_1) +   \dots + \widehat{\va_i
  (C_i)}+
  \dots +\va_k(C_k))(x_{e})=0$
  \item [(ii)] $\nabla_{xx}^2 (C_i + \va_1 (C_1) + \dots+
  \widehat{\va_i (C_i)}+ \dots+ \va_k(C_k))(x_e)_{|_{\widetilde{W}
  \times \widetilde{W}}}$ is positive or
  negative definite, where $$\widetilde{W} := \bigcap\limits_
  {\tiny \begin{array}{l} j=1\\j \not = i\end{array}}^k
   \ker (d \va_j (C_j)) (x_e),$$
\end{itemize}
then $x_e$ is nonlinear stable.
\end{theorem}

The aim of our paper is to prove the equivalence of these three
methods. This shows that when $x_e$ is an equilibrium point for
(\ref{11}) and we choose $C_1, \dots, C_k$ as a set of constants
of motion, if we conclude stability of $x_e$ with one of the
methods, then the other two will also give stability of $x_e$.
Thus we can choose the most convenient method from the
computational point of view. Since computations can become
cumbersome in some examples it is important to know that if we
cannot conclude stability of $x_e$ using the set $C_1, \dots, C_k$
of constants of motion with one of the methods, then we cannot
conclude stability of $x_e$ by applying the other two methods
using the same set $C_1, \dots, C_k$ of constants of motion.

\section{Equivalence of the three methods}

In order for the paper to be self-contained we will start by
proving Arnold's result on stability since in his original paper
\cite{1} the proof was omitted. In order to do this we need the
following preliminary results which will play a crucial role in
all that follows.

We will begin by establishing the notations and conventions to be
used throughout this paper. A vector $x \in \mathbb{R}^n$ will be
considered as a column vector or a $n \times 1$ matrix. Its
transpose will be a row vector or a $1 \times n$ matrix.

Let $f:\mathbb{R}^n \to \mathbb{R}$ be a $C^1$ real valued
function. The gradient of $f$ at a point $x \in \R^n$ is defined
as the column vector
\begin{equation*}
    \nabla f(x) = \left[\begin{array}{c} \ds\frac{\pa f}{\pa x_1}
    (x)\\ \vdots \\ \ds\frac{\pa f}{\pa x_n} (x) \end{array}\right].
\end{equation*}
If $f :\R^n \to \R^m$ is a vector valued map, then it will be
represented as a column vector of its component functions $f_1,
\dots, f_m$, namely
\begin{equation*}
    f(x) = \left[\begin{array}{c} f_1 (x)\\ \vdots\\f_m (x)
    \end{array}\right].
\end{equation*}

If $f\in C^1 (\R^n, \R^m)$, then we introduce the notation
\begin{equation*}
    \nabla f(x) : = \left[ \nabla f_1 (x) \dots \nabla f_m (x)
    \right],
\end{equation*}
where $\nabla f(x)$ is a $n \times m$ matrix which has  as columns
the gradient vectors $\nabla f_1 (x), \dots, \nabla f_m (x)$. Note
that the transpose matrix $\nabla f(x)^T$ is the  Jacobian matrix
of $f$ at the point $x \in U$.

Let $f: \R^{n+k}\to \R$ be a $C^2$ real valued function and $(x,y)
\in \R^n \times \R^k$. We will use the following notations,
\begin{equation*}
    \begin{array}{ll} \nabla_x f(x, y) =\left[\begin{array}{c}\ds\frac{\pa f
    (x,y)}{\pa x_1}\\ \vdots\\ \ds\frac{\pa f(x, y)}{\pa
    x_n}\end{array}\right], &  \nabla_y f(x, y) =\left[\begin{array}{c}\ds\frac{\pa f
    (x,y)}{\pa y_1}\\ \vdots\\ \ds\frac{\pa f(x, y)}{\pa
    y_k}\end{array}\right],\end{array}
\end{equation*}

\begin{equation*}
    \begin{array}{ll}
    \nabla_{xx}^2 f(x, y) = \left[\ds\frac{\pa^2 f(x,y)}{\pa x_i \pa
    x_j}\right], & \nabla_{xy}^2 f(x, y) = \left[\ds\frac{\pa^2 f(x,y)}{\pa x_i \pa
    y_j}\right],\\\\  \nabla_{yy}^2 f(x, y) = \left[\ds\frac{\pa^2 f(x,y)}{\pa y_i \pa
    y_j}\right].\end{array}
\end{equation*}

For the proof of Theorem \ref{tA} we will need the following
result that can be found in references \cite{2} and \cite{6} .

\begin{proposition}\label{p21}
Let $P$ be a symmetric $n \times n$ matrix and $Q$ a positive
semidefinite symmetric $n \times n$ matrix. We assume that
\begin{equation*}
    x^T Px >0,
\end{equation*} for all $x \in \R^n$, $x \not = 0$ satisfying $x^T
Qx=0$. Then there exists a scalar $\al \in \R$ such that
\begin{equation*}
    P+\al Q >0.
\end{equation*}
\end{proposition}

\begin{proof}
We will prove by contradiction. Then  for every integer $k$, there
exists a vector $x_k \in \R^n$ with $\|x_k\|=1$ such that:
\begin{equation}\label{22}
x_k^T P x_k + k x_k^T Qx_k \leq 0.
\end{equation}

The sequence $(x_k)$ is bounded and consequently it has a
subsequence, that we will denote also by $(x_k)$, converging to a
vector $\ol x \in \R^n$ with $\|\ol x\|=1$. Taking the limit in
\eqref{22} we obtain
\begin{equation}\label{23}
\ol x^T P\ol x + \lim\limits_{k \to \infty} (kx_k^T Q x_k) \leq 0.
\end{equation}
Since \begin{equation*}
    x^T_k Qx_k \geq 0,
\end{equation*} the inequality \eqref{23} implies that $(x_k^TQx_k)$
converges to zero and hence $\ol x^T Q \ol x=0$.

It follows from the hypothesis  that $\ol x^T P \ol x >0$ and this
contradicts \eqref{23}.
\end{proof}

\medskip

Let $x_e$ be an equilibrium point for the dynamic \eqref{11} and
let $C_1, \dots, C_k \in C^1 (\R^m, \R)$ be a set of constants of
notion for the dynamic \eqref{11}. We define the following
quadratic form,
\begin{equation}\label{24}
x^T Q_{i} x: = \sum_{\tiny \begin{array}{l} j=1\\j
    \not = i\end{array}}^k x^T \nabla C_j (x_e)(\nabla C_j (x_e))^T x.
\end{equation}

We have the following characterization for the vector subspace $W$
defined in Theorem \ref{tA}.

\begin{lem} \label{l25}
\begin{eqnarray*} x^TQ_{i}x =0 \Leftrightarrow x \in W=
\bigcap\limits_{\tiny \begin{array}{l}
  j=1\\j \not = i\end{array}}^k \ker d C_j.\end{eqnarray*}
\end{lem}

\begin{proof}
\begin{eqnarray*}
  x^TQ_{i}x =0 &\Leftrightarrow&  \sum_{\tiny \begin{array}{l} j=1\\j
    \not = i\end{array}}^k x^T \nabla C_j (x_e)(\nabla C_j (x_e))
  ^T x=0\\&\Leftrightarrow& \sum_{\tiny \begin{array}{l} j=1\\j
    \not = i\end{array}}^k((\nabla C_j (x_e))^T x)
  ((\nabla C_j (x_e))^T x)=0\\&\Leftrightarrow&  \sum_{\tiny \begin{array}{l} j=1\\j
    \not = i\end{array}}^k((\nabla C_j
  (x_e))^T x)^2=0\\&\Leftrightarrow&  (\nabla C_j (x_e))^T x=0, \quad \forall
  \ j=\ol{1,k}, j\not =i \\&\Leftrightarrow& x \in W.
\end{eqnarray*}

\end{proof}

\medskip

{\bf Proof of Theorem \ref{tA}}. Let $L_{i, \al_i}\in C^2 (\R^n
\times \R^{k-1}, \R)$ be the function defined by
\begin{eqnarray*}
    L_{i, \al_i}  (x, \lambda_1, \dots, \widehat{\lambda_i}, \dots, \lambda_k)&: =& F_i (x, \lambda_1, \dots, \widehat{\lambda_i}, \dots,
    \lambda_k)\\&+&\ds\frac{\al_i}{2}\sum_{\tiny \begin{array}{l} j=1\\j
    \not = i\end{array}}^k \left[C_j(x) - C_j (x_e)\right]^2,
\end{eqnarray*} where $\al_i \in \R$ will be determined later.

A simple computation shows us that
\begin{eqnarray*}
    \nabla_x L_{i, \al_i} (x, \lambda_1, \dots, \widehat{\lambda_i}, \dots,
    \lambda_k)&=&\nabla_xF_i (x, \lambda_1, \dots, \widehat{\lambda_i}, \dots,
    \lambda_k)\\&+&\al_i \sum_{\tiny \begin{array}{l} j=1\\j
    \not = i\end{array}}^k (C_j(x)-C_j(x_e))\nabla C_j (x)
\end{eqnarray*} and
\begin{eqnarray*}
    \nabla_{xx} L_{i, \al_i} (x_e , \lambda_1, \dots, \widehat{\lambda_i}, \dots,
    \lambda_k)&=&\nabla_{xx} F_i (x_e, \lambda_1, \dots, \widehat{\lambda_i}, \dots,
    \lambda_k)\\&+&\al_i \sum_{\tiny \begin{array}{l} j=1\\j
    \not = i\end{array}}^k  \nabla_x C_j(x_e) (\nabla
    C_j(x_e))^T\\&=& P_i +\al_i Q_i,
\end{eqnarray*} where $Q_i=\sum\limits_{\tiny \begin{array}{l} j=1\\j
    \not = i\end{array}}^k \nabla_x C_j(x_e)(\nabla C_j(x_e))^T$ is
    the $n \times n$ symmetric matrix  defined by \eqref{24}.

The hypothesis (i) implies that $\nabla_x L_{i, \al_i}(x_e,
\lambda_1^*, \dots, \widehat{\lambda^*_i}, \dots, \lambda^*_k)=0$.
As a consequence of the hypothesis (ii) and  Proposition \ref{p21}
we can find $\al_i^* \in \R$ such that $P_i +\al_i^* Q_i >0$ and
thus $L_{i,\al_i^*}(x)> 0$ for $x\neq x_{e}$ in a small
neighborhood of the equilibrium point $x_{e}$.

Let us define now the function $V_{i, \al_i^*} \in C^2 (\R^n, \R)$
by the following relation,
\begin{equation*}
    V_{i, \al_i^*}(x) = L_{i, \al_i^*} (x, \lambda_1^*, \dots,
    \widehat{\lambda_i^*}, \dots, \lambda_k^*) - L_{i, \al_i^*} (x_e,
    \lambda_1^*, \dots, \widehat{\lambda_i^*}, \dots, \lambda_k^*).
\end{equation*} It is easy to see that $V_{i, \al_i^*}$ is a
Lyapunov function and consequently via Lyapunov's theorem the
equilibrium state $x_e$ is nonlinear stable.\quad $\blacksquare$

\medskip

The proofs of Theorem \ref{tB} and Theorem \ref{tC} can be found
in the original papers \cite{3}  and \cite{5}. They are also based
on finding a corresponding Lyapunov function.

Now we will prove the main result of this paper.

\begin{theorem}\label{t25}
 Let $C_1, \dots, C_k \in C^2 (\R^n, \R)$ be a
 set of constants of motion for the dynamic (\ref{11}).
 Then the following statements are equivalent:
\begin{itemize}
  \item [(a)] hypotheses of Theorem \ref{tA} hold;
  \item [(b)] hypotheses of Theorem \ref{tB} hold;
  \item [(c)] hypotheses of Theorem \ref{tC} hold.
\end{itemize}
Each of the above statements implies nonlinear stability.
\end{theorem}

\begin{proof}
"$(a)\Rightarrow (b)$" Assume that the hypotheses of Theorem
\ref{tA} hold. Consider the following functions $\va_j: \R \to
\R$, $\va_j (t) = - \lambda_j^* t + \ds\frac{\al_i}{2} (t- C_j
(x_e))^2$, for $j \not = i$ and $\al_j \in \R$ arbitrary for the
moment, and $\lambda_j^*$ given in Theorem \ref{tA}. As in  the
proof of Theorem \ref{tA}, the conditions (i) and (ii) of Theorem
\ref{tA} imply the conditions (i) and (ii) of Theorem \ref{tB} for
a certain $\al_j^*$ given by Proposition \ref{p21}.

"$(b)\Rightarrow (c)$" This is obvious since positive or negative
definiteness on the whole space implies positive or negative
definiteness on the subspace $\widetilde{W}$.

"$(c)\Rightarrow (a)$". Assume that the hypotheses of Theorem
\ref{tC} hold. Let $\lambda_j^* = - \va_j^\p (C_j (x_e))$ for $j
\not = i$. It is obvious that condition (i) of Theorem \ref{tC}
implies condition (i) of Theorem \ref{tA}. Also because some of
$\lambda_j^*$'s might be zero we have the inclusion $W \subseteq
\widetilde{W}$. Then
\begin{eqnarray*}
  && z^T \left[ \nabla_{xx}^2 (C_i +\va_1 (C_1) + \dots + \widehat
  {\va_i (C_i)} +\dots
  + \va_k (C_k))(x_e)\right] y  \\
  &=& z^T \nabla_{xx}^2 C_i(x_{e})y+\sum_{\tiny \begin{array}{l}
   l=1\\l \not = i \end{array}}
  ^k z^T \left( \va^\p_l (C_{l}(x_e))
  \left[\ds\frac{\pa^2 C_l(x_{e})}{\pa x_i \pa x_j}\right]
  \right)y \\ &+& \sum_{\tiny \begin{array}{l} l=1\\l \not = i \end{array}}^k
  \sum_{s,p=1}^n \va^{\p \p}_l (C_l (x_e)) z_s y_p \ds\frac{\pa C_l(x_{e})}
  {\pa x^s} \ds\frac{\pa C_l(x_{e})}{\pa x^p} \\
   &=& z^T \nabla_{xx}^2 F_i (x_e, \lambda_1^*, \dots,
   \widehat{\lambda_i^*}, \dots, \lambda_k^*)y+\sum_{\tiny \begin{array}{l}
   l=1\\l \not = i \end{array}}^k
  \sum_{s,p=1}^n \va^{\p \p}_l (C_l (x_e)) z_s y_p \ds\frac{\pa C_l(x_{e})}
  {\pa x^s} \ds\frac{\pa C_l(x_{e})}{\pa
  x^p},
\end{eqnarray*}
 for any $z, y \in \R^n$.

 If we take $z, y \in W$ the second summand will be zero and
 consequently condition (ii) of Theorem \ref{tC} implies condition (ii) of
 Theorem \ref{tA}.
\end{proof}

\medskip

In all of the three methods the stability is decided when a
certain matrix is positive or negative definite. Consequently,
Arnold's method seems to be the most economical since it requires
definiteness of a smaller matrix than the other two methods.

Next we will discus the situation in which condition (i) of
Theorem \ref {tA} is not satisfied. Or equivalently, when the
vectors $\nabla C_{i}(x_{e})$, $i\in \overline{1,k}$ are linear
independent. Consequently, in a small neighborhood $U_{x_{e}}$ of
$x_{e}$ they generate an integrable distribution whose leaves are
the level sets of the map $F:=(C_1, \dots ,C_k):\mathbb{R}^n
\rightarrow \mathbb{R}^k$. Eventually after shrinking $U_{x_{e}}$
all the points in $U_{x_{e}}$ are regular points for $F$. There
exists a diffeomorphism $\phi:U_{x_{e}}\rightarrow
(F^{-1}(F(x_{e}))\cap U_{x_{e}} )\times V_{F(x_{e})}$, where
$V_{F(x_{e})}$ is a small neighborhood of $F(x_{e})$ in
$\mathbb{R}^k$. Because $(C_1, \dots ,C_k)$ are constants of
motion for the dynamic (\ref{11}) we obtain $\phi_{\ast}f=(Y,0)$,
where $Y\in \mathfrak{X}\left( F^{-1}(F(x_{e}))\cap
U_{x_{e}}\right)$. If $(y,z)$ are coordinates induced by $\phi$ on
$(F^{-1}(F(x_{e}))\cap U_{x_{e}} )\times V_{F(x_{e})}$ from a set
of coordinates around $x_e$ then the equations of motion
corresponding to the vector field $\phi_{\ast}f$ are

\begin{equation}\label{bifu}
\begin{array}{l} \dot y=Y(y,z)\\\\ \dot z=0.
\end{array}
\end{equation}
Moreover, $\phi(x_{e})=(y_{e},0)$ and $y_{e}$ is an equilibrium
point for the dynamic generated by the vector field $Y$. The above
system can be regarded as a bifurcation problem with $z\in
V_{F(x_{e})}$ the bifurcation parameter. We have the following
result.

\begin{theorem}
If $(C_1, \dots ,C_k)$ are constants of motion for the dynamic
(\ref{11}) which are linear independent at the equilibrium point
$x_{e}$, then $x_{e}$ is stable for the dynamic (\ref{11}) if the
equilibrium point $y_{e}$ is stable for the dynamic generated by
the vector field $Y$ and $(y_{e},0)$ is not a bifurcation point
for (\ref{bifu}).
\end{theorem}

\medskip

This result was used in \cite{bipu} for the stability problem of
Ishii's equation. Given the conditions of the above theorem it is
enough to study the stability of a dynamical system that has fewer
variables. Nevertheless, the problem is not free of difficulties
since one has to find a set of adapted coordinates around $x_{e}$
for the local fibration generated by the map $F$.

\section{Examples}

\subsection{The free rigid body}

Theorem \ref{t25} asserts that if  stability is obtained with one
of the methods, then it can be obtained with the other two as
well. Indeed, let us consider the Euler momentum equations:

\begin{equation*}
    \left\{\begin{array}{l} \dot m_1 =
    \left(\ds\frac{1}{I_3}-\ds\frac{1}{I_2}\right) m_2 m_3\\\\ \dot
    m_2 = \left(\ds\frac{1}{I_1}-\ds\frac{1}{I_3}\right) m_1m_3\\\\
    m_3 = \left(\ds\frac{1}{I_2}-\ds\frac{1}{I_1}\right) m_1m_2
    \end{array}\right.
\end{equation*}
where $I_1 >I_2 >I_3 >0$. Then $x_e = (M, 0,0)$ is an equilibrium
point and $C_1 (m_1, m_2, m_3) = \ds\frac{1}{2}
\left(\ds\frac{m_1^2}{I_1}+\ds\frac{m_2^2}{I_2}+\ds\frac{m_3^2}{I_3}\right)$,
and $C_2 (m_1, m_2, m_3) = \ds\frac{1}{2}(m_1^2+m_2^2+m_3^2)$ are
two constants of motion.

We study the stability of $x_e= (M, 0,0)$, $M \not = 0$ by using
Arnold's method. Let $F_1 = C_1 - \lambda C_2$, then $\nabla F_1
(x_e) =0$ iff $\lambda = \ds\frac{1}{I_1}$. Also
\begin{equation*}
    \nabla_{xx}^2 F_1 \left(x_e, \ds\frac{1}{I_1}\right) =
    \left[\begin{array}{ccc}0&0&0\\0&
    \ds\frac{1}{I_2}-\ds\frac{1}{I_1}&0\\0&0&\ds\frac{1}
    {I_3}-\ds\frac{1}{I_1}\end{array}\right]
\end{equation*} and $W = Span \left( \left[\begin{array}{c}
0\\1\\0\end{array}\right], \left[\begin{array}{c}
0\\0\\1\end{array}\right]\right)$. It is easy to see that
$$\nabla_{xx}^2 F_1 \left(x_e , \ds\frac{1}{I_1}\right)_{|_{W\times
W}}>0$$. This shows that $x_e = (M, 0,0)$, $M \not = 0$ is
nonlinear stable.

Next, we will prove the same stability result using the other two
methods. We begin with the Energy-Casimir method. Let $H_\va = C_1
+ \va (C_2)$. The first variation is given by
\begin{eqnarray*}
   \delta H_\va&=& \ds\frac{m_1}{I_1}\de m_1 + \ds\frac{m_2}{I_2}
   \de m_2 +\ds\frac{m_3}{I_2}\de m_3
   + \va^\p (m_1 \de m_1 + m_2 \de m_2 + m_3 \de m_3).
\end{eqnarray*}

Then $\de H_\va (M,0,0) =0$ is equivalent with $\va^\p
\left(\ds\frac{1}{2}M^2\right) = -\ds\frac{1}{I_1}$. Also
\begin{eqnarray*}
   \delta^2 H_\va (M, 0,0)&=& \left(\ds\frac{1}{I_2}
   -\ds\frac{1}{I_1}\right) (\de m_2)^2 +
   \left(\ds\frac{1}{I_3} -\ds\frac{1}{I_1}\right)(\de m_3)^2  \\
   &+& \va^{\p \p} \left( \ds\frac{1}{2} M^2\right) M^2 (\de m_1)^2
\end{eqnarray*} is positive definite iff $\va^{\p \p}
\left(\ds\frac{1}{2}M^2\right) >0$.

We can take $\va (t) = \left(t - \ds\frac{1}{2}M^2\right)^2 -
\ds\frac{1}{I_1}t$ and  conclude that $x_e = (M, 0,0)$, $M \not =
0$ is nonlinear stable.

For Ortega-Ratiu's method we can take the same constant of motion
used for applying Arnold's method, i.e. $F_1 =C_1-\ds\frac{1}{I_1}
C_2$.

\subsection{Lorenz five component model}

We will show in this example that if the stability of an
equilibrium point cannot be decided with one of the three methods
then it cannot be decided with the other two either. This is what
Theorem \ref{t25} is predicting. It simplifies enormously the
computations in the sense that if we do the computations using one
of the methods and obtain that the stability cannot be decided,
then it is useless to do the computations using the other two
methods and the same set of constants of motion.

To illustrate this, we will take the example of Lorenz five
component model. The equations are

 \begin{equation*}
    \left\{\begin{array}{l} \dot x_1 = -x_2 x_3 + bx_2x_5 \\ \dot
    x_2 = x_1 x_3 - bx_1 x_5\\ \dot x_3 = -x_1x_2\\ \dot x_4 =
    -\ds\frac{x_5}{\varepsilon}\\ \dot x_5 =
    \ds\frac{x_4}{\varepsilon} + bx_1 x_5\end{array}\right.
 \end{equation*} where $b, \varepsilon \in \R^{\ast}$,
 $x_e = (0, 0, M, 0, 0)$, $M \not = 0$ is an  equilibrium point and
$C_1 (x_1, \dots, x_5)= \ds\frac{1}{2}(x_1^2 + 2 x_2^2 + x_3^2 +
x_4^2 + x_5^2)$, and $C_2 = \ds\frac{1}{2}(x_1^2 + x_2^2)$ are
constants of motion.

We try to apply Arnold's method. Take $F_1 = C_1 - \lambda C_2$.
Then $\nabla F_1 (x_e) = 0$ is impossible for any $\lambda \in
\R$. We have another possibility for choosing a constant of
motion. Let $F_2 = C_2 - \lambda C_1$. Then $\nabla F_2 (x_e) =0$
iff $\lambda = 0$. Also
\begin{equation*}
    \nabla_{xx}^2 F_2 \left(x_e, 0\right) =
    \left[\begin{array}{ccccc}1&0&0&0&0\\0&
    1&0&0&0\\0&0&0&0&0\\0&0&0&0&0\\0&0&0&0&0 \end{array}\right]
\end{equation*} and $$W = Span \left( \left[\begin{array}{c}
1\\0\\0\\0\\0\end{array}\right], \left[\begin{array}{c}
0\\1\\0\\0\\0\end{array}\right], \left[\begin{array}{c}
0\\0\\0\\1\\0\end{array}\right], \left[\begin{array}{c}
0\\0\\0\\0\\1\end{array}\right]\right).$$

It is easy to see that $\nabla_{xx}^2 F_2 (x_e, 0)_{|_{W\times
W}}$ is not definite.

Now we try to apply the Energy-Casimir method. Let $H^1_\va = C_1
+ \va (C_2)$. Then $\de H_\va^1 (x_e) = 0$ is impossible for any
$\va \in C^2 (\R, \R)$. We take the other possibility, namely
$H^2_\va = C_2 + \va (C_1)$. Then we have
\begin{equation*}
    \de H^2_\va = x_1 \de x_1 + x_2 \de x_2 + \va^\p (x_1 \de x_1 +
    2 x_2 \de x_2 + x_3 \de x_3 + x_4 \de x_4 + x_5 \de x_5).
\end{equation*}
Consequently $\de H^2_\va (x_e) = 0$ iff $\va^\p
\left(\ds\frac{1}{2}M^2\right) =0$. Also
\begin{equation*}
    \de^2H^2_\va = (\de x_1)^2 + (\de x_2)^2 + \va^{\p \p}
    \left(\ds\frac{1}{2}M^2\right) (\de x_3)^2
\end{equation*} which is not definite.

Finally we will try to apply Ortega-Ratiu's method. Let $F= C_2 +
\va (C_1)$. We have that $\de F(x_e) =0$ iff $\va^\p
\left(\ds\frac{1}{2}M^2\right) =0$ and then $\widetilde{W} =
\mathbb{R}^5$. Also
\begin{equation*}
    \de^2 F (x_e) =
    \left[\begin{array}{ccccc}1&0&0&0&0\\0&
    1&0&0&0\\0&0&\va^{\p \p} \left(\ds\frac{1}{2}M^2\right)
    M^2&0&0\\0&0&0&0&0\\0&0&0&0&0 \end{array}\right]
\end{equation*} and consequently $\de^2 F(x_e)_
{|_{\widetilde{W}\times\widetilde{W}}}$ is not definite for any
choice of  $\va \in C^2 (\R, \R)$.

\noindent {\sc P. Birtea} \\
Departamentul de Matematic\u a, Universitatea de Vest,
RO--1900 Timi\c soara, Romania.\\
D\'epartement de Math\'ematiques de Besan\c{c}on, Universit\'e de
Franche-Comt\'e, UFR des Sciences et Techniques, 16 route de Gray,
F--25030 Besan\c{c}on c\'edex,
France.\\
Email: {\sf birtea@math.uvt.ro}
\medskip

\noindent {\sc M. Puta} \\
Departamentul de Matematic\u a, Universitatea de Vest,
RO--1900 Timi\c soara, Romania.\\
Email: {\sf puta@math.uvt.ro}
\end{document}